\input amstex 
\documentstyle{amsppt}
\input bull-ppt
\keyedby{bull436e/lic}

\define\SPAN{\text{span}}
\define\CH{\text{ch}}
\define\OV{\overline{V}}

\define\OW{\overline{W}}
\define\Spec{\text{Spec}}
\define\rank{\text{rank}}
\define\Int{\text{Int}}

\topmatter
\cvol{30}
\cvolyear{1994}
\cmonth{January}
\cyear{1994}
\cvolno{1}
\cpgs{62-69}
\ratitle
\title Genera of Algebraic Varieties\\and Counting of 
Lattice Points \endtitle
\author Sylvain E. Cappell and Julius L. Shaneson \endauthor
\address   Courant Institute,  New York University, New 
York, New York 10012\endaddress
\address Department of Mathematics, University of 
Pennsylvania, Philadelphia,
Pennsylvania 19104\endaddress
\ml jshaneso\@mail.sas.upenn.edu\endml
\settabs 2 \columns
\date October 16, 1992 and, in revised form, January 4, 
1993\enddate
\thanks Both authors were partially supported by NSF 
grants\endthanks
\abstract This paper announces results on the behavior of 
some important
algebraic and topological invariants --- Euler 
characteristic,
arithmetic genus, and their intersection homology 
analogues; the signature,
etc. --- and their associated characteristic classes, 
under morphisms
of projective algebraic varieties. The formulas obtained 
relate global
invariants to singularities of general complex algebraic 
(or analytic)
maps. These results, new even for complex manifolds, are 
applied to obtain
a version of Grothendieck-Riemann-Roch, a calculation of 
Todd classes
of toric varieties, and an explicit formula for the number 
of
integral points in a polytope in Euclidean space with 
integral vertices.\endabstract
\endtopmatter

\document

Consider first the behavior of the classical 
Euler-Poincare characteristic
$$
e (X) =\sum _ i (-1)^i \rank \, H^i(X)
$$
under a (surjective) projective morphism $f\colon\, X \to Y$
of projective (possibly singular) algebraic varieties. 
Such a morphism
can be stratified with subvarieties as strata. In 
particular, there is
a filtration of $Y$ by closed subvarieties, underlying a 
Whitney
stratification,
$$
\phi\subset Y _0 \subset \cdots \subset Y _s = Y
$$
of strictly increasing dimension, such that $Y _i - 
Y_{i-1}$ is a union
of smooth
 manifolds of the same dimension and such that the 
restriction
of $f$ to $f^{-1}(Y _i - Y_{i-1})$ is a locally trivial 
map of Whitney
stratified spaces. (In the results it will suffice to have 
$\dim \, Y_i < \dim\, Y_{i+1}$ and $\dim \, f^{-1}(x)$ 
constant
over \RM{``}strata\RM{''} $Y_i - Y_{i-1}$.)

We recall the definition of the normal cone $C_Z W$ of an 
irreducible
subvariety $Z$ of a variety $W$:
$$
C_ZW =\Spec \left(\bigoplus {\Cal I} ^{n}/{\Cal I} ^{n+
1}\right)\ ,
$$
${\Cal I}$ the sheaf of ideals defining $Z$. Let 
$P(C_Z\oplus 1)$ be
its projective completion, and let $P_{Z,W}$ be the 
general fiber
of the canonical projective morphism [F1]
$$
P(C_Z\oplus 1)\to Z \ .
$$
For example, if $Z$ is a smooth subvariety of $W$ of 
complex codimension 
$d$, $P_{Z,W}$ is the complex projective space $CP^d$ of 
dimension $d$.

Let ${\Cal V}$ be the set of components of strata of $Y$. 
For $V\in {\Cal V}$
define $\hat{e}(\OV)$ inductively by the formula
$$
\hat{e}(\OV) = e (\OV) -\sum_{W<V} \hat{e}
(\OW)\, e (P_{W,V})\ ,
$$
where the sum is over all $W\in {\Cal V}$ with 
$W\subset\OV-V$.
Let $P_{V,f}$ be the general fiber of the composite
$$
P(C_{f^{-1}V} \oplus 1 )\longrightarrow 
f^{-1}V\longrightarrow V\ .
$$

\proclaim{Theorem 1}
Let $f\colon\, X\to Y$ be a surjective morphism of 
projective complex
algebraic varieties, and let ${\Cal V}$ be the components 
of the 
\RM{``}strata\RM{''}.
Assume that $\pi _1V=0$ for $V\in {\Cal V}$. Let $F$ be 
the general fiber
of $f$. Then
$$
e (X) = e (Y)\, e (F)+\sum_{{\Cal V}_0} \widehat{e}(\OV)
\left[ e (P_{V,f}) - e (F)\, e(P_{V,Y})\right]\ ,
\tag"$(\ast)$"
$$
$$
{\Cal V}_0 \text{ the subset of } {\Cal V} \text{ with }
\dim \, V< \dim \, Y \ .
$$
\endproclaim

\ex{Example} Let $X$ be obtained from $Y^n$ by blowing
up a point $y$. Let $D=P(C_{\{y\}})=f^{-1}y\subset X$ be 
the exceptional
set. Then $(*)$ becomes
$$
e (X) = e (Y) + 2e (D) -e (P(C_{\{y\}}\oplus 1)\ .
$$
If $y$ is a smooth point, $D$ and $P(C_{\{y\}}\oplus 1)$ 
are complex
projective spaces $CP^{n-1}$ and $CP^n$, and this is the 
well-known result
$$
e (X)= e (Y) + n -1 \ .
$$
\endex

\dfn{Definition} Any invariant satisfying $(*)$ under the
hypotheses of Theorem~1 will be said to have the 
stratified multiplicative
property (SMP).
\enddfn

Without the $\pi _1$ hypothesis the SMP and the present 
results must be
phrased in terms of coefficients in local systems. 
Theorem~1 for this
case essentially includes as a corollary results like the 
generalization
of Riemann-Hurwitz given in \cite{DK; I; K, (III, 32)}.

As another example consider the signature $\sigma(X)$, 
defined in
[GM2] as the signature (number of positives minus number 
of negatives
in a diagonalization) of the intersection form in the 
middle-dimensional intersection homology 
$\operatorname{IH}^{\overline m}_n(X)$ with
middle perversity, $X$ projective of complex dimension $n$.
We show the signature also has the SMP. \ For blowing up a 
point
this becomes
$$
\sigma (X) =\sigma(Y) -\sigma (P(C_{\{y\}} \oplus 1)\ .
$$

In [CS1] a different formula was given for the behavior of 
the signature
(and the $L$-classes) under any stratified map. In the 
case of blowing
up a point the formula of [CS1] gives
$$
\sigma (X)=\sigma (Y) +\sigma(E _y)\ ,
$$
where
$$
E _y = f^{-1}(N) \, /\, f^{-1}(L)\ ,
$$
$N$ a (piecewise linear) neighborhood of $y\in V$ and 
$L=\partial N$ its
frontier. There are two key differences between this 
formula and the above.
The first is that the topological completion $E_y$ of 
$f^{-1}(\Int\, N)$ 
usually is not an algebraic variety. The second is that, 
unlike $E _y$,
it can and usually does happen that no open neighborhood 
of $f^{-1}(y)$
in $P_{V,f}$ (which is actually a bundle over $D$) embeds 
topologically
in $X$, extending the inclusion of $f^{-1}(y)\subset X$.
Thus the behavior of genera under {\it algebraic} maps can 
be
determined without knowing the local topological structure 
around
strata precisely.

We now consider two extensions of Hirzebruch genera of 
smooth
varieties to general varieties. Let $Ih^{p,q}(X)$ be the 
Hodge
numbers of Saito's pure Hodge decompositon [S] on 
$\operatorname{IH}^{\overline m}_\ast (X;{\Bbb C})$, and 
let 
$h^{p,q}_i(X)$ be the Hodge numbers of
Deligne's mixed Hodge structure on $H^i(X;{\Bbb C})$. Let 
$y$
be a variable, and define genera and intersection genera by
$$
\chi _y (X) = \sum _p\left[ \sum_{i,q} (-1)^{i-p}\, 
h^{p,q}_i(X)\right]y^p\ 
$$
and
$$
I\chi _y (X) = \sum _p\left[ \sum_{q} (-1)^{q}\, 
Ih^{p,q}(X)\right]y^p\ .
$$
Thus $\chi_{_{-1}}=e$, $I\chi _{_1} =\sigma$, and $I\chi 
_{_0}$ and $\chi _0$
are two possible extensions to singular varieties of the
arithmetic genus.

\proclaim{Theorem 2}
The genera $\chi _y$ and $I \chi _y$ have the SMP.
\endproclaim

We discuss briefly the proof of Theorem~2 for the 
intersection genera.
The strategy is to argue inductively, blowing up closures 
of strata
starting with a lowest dimensional component, relying on 
the following
result on general blowing up:

\proclaim{Theorem 3}
Let $X$ be a projective complex algebraic variety and $Z$ 
a closed
subvariety. Let $\widetilde{X} = Bl _Z X$ be the blowup of 
$X$
along $Z$, with exceptional divisor $E=P(C _ZX)$. Then
$$
I\chi _y(\widetilde{X}) 
= I \chi _y(X) + (1-y) \chi _y (E)-I\chi _y (P(C _Z 
X\oplus 1))\ .
$$
\endproclaim

This is proved by considering a diagram:
$$
\CD
\widetilde{X} @>\tilde\psi>>
\widetilde{X} \cup _E P(C _E \widetilde{X} \oplus 1)\\
@VVV @VVV\\ 
X @>\psi>>
\widetilde{X} \cup_ E P(C _Z X\oplus 1)
\endCD
$$
The vertical maps are the projection from the blowups 
along $Z\subset X$
and along the canonical embedding
$$
Z\subset C _Z X \subset P(C _Z X\oplus 1) \subset 
\widetilde{X} \cup_E
P(C_Z X\oplus 1)\ ;
$$
note that
$$
P(C _E\widetilde{X} \oplus 1) = Bl_Z P(C _Z X\oplus 1)\ .
$$
The horizontal arrows are the specialization maps from the 
general
to the special fiber in the deformation to the normal cone 
\cite{H, V,
GM1}. For the top map note that $Bl 
_E\widetilde{X}=\widetilde{X}$.
The main idea is to compare the monodromy weight filtrations
and associated spectral sequences for $R\psi _\ast 
{\Bbb I} {\Bbb C}^\bullet (X)$ and $R\tilde\psi _\ast 
{\Bbb I}{\Bbb C}^\bullet (\widetilde{X})$, taking 
advantage of
the result of [S], strengthening [BBD], that the vertical 
maps
``induce'' injections in intersection cohomology onto 
summands,
respecting pure Hodge structures.

For each genus or intersection genus we define corresponding
characteristic homology class $T _y$ or $IT _y$. For 
example,
$T_{-1}$ is the total MacPherson Chern class, $T _0 = Td$
(the image in homology of) the Todd class of [BFM], 
appearing
in the generalized GRR theorem (see also [F1]), and $IT_1$
is the Goresky-MacPerhson $L$-class. We show these also have
the SMP; i.e., they satisfy the appropriate version of 
$(*)$.
Combining this with Grothendieck-Riemann-Roch yields, for 
example:

\proclaim{Corollary}
Let $f\colon\, X\to Y$ be as in Theorem~\RM1. Let $\alpha$ 
be a
locally free sheaf on $X$, and let $f _\ast \alpha$ be a 
locally
free sheaf on $Y$ with
$$
f _\ast (\alpha\otimes {\Cal O} _X) =f _\ast \alpha \otimes 
{\Cal O} _Y \ .
$$
\RM(This always exists for $Y$ smooth.\RM) Then
$$
\aligned
&f _\ast Td _X -\chi (F;{\Cal O}_F)\, \CH^{-1}(f _\ast 
\alpha)\cap f _\ast
(\CH (\alpha)\cap Td _X) \\
&\qquad = \sum_{{\Cal V}_0} i _\ast \widehat{T} d _{\OV} 
\left[\chi
(P_{V,f}; {\Cal O}_{P_{V,f}})-\chi (F;{\Cal O}_F)\, 
\chi(P_{V,Y};
{\Cal O}_{P_{V,Y}})\right] \ . 
\endaligned
$$
\endproclaim

Finally, we consider toric varieties and counting of 
lattice points.
Toric varieties arise naturally in algebraic and 
symplectic geometry,
representation theory, number theory, and combinatorics. A 
toric
variety $X^n$ is an irreducible normal variety on which 
the complex
torus $({\Bbb C}-\{0\})^n$ acts with an open orbit
\cite{O, D}. Equivalently, it is given locally by systems 
of monomial
equations. If $X$ is projective, the quotient of $X$ by 
the action
of the real torus $(S^1)^n$ can be naturally identified 
(via the
``moment map'' if $X$ happens to be a symplectic manifold) 
with a
convex polytope in ${\Bbb R}^n$ whose vertices have 
integral coordinates,
and questions about the polytope can often be translated 
into questions
about the variety and vice~versa.

Given the polytope $P$ with vertices lattice points, there 
is a classical
recipe for constructing the corresponding toric variety. 
This is often
done via the fan $\Sigma$ dual to the polytope; in fact, 
every toric 
variety has such a description uniquely in terms of a fan 
\cite{O, 1.2}.
A complete fan is a decomposition of ${\Bbb R}^n$ into 
cones, each
spanned by its ``rays'', i.e., half lines from the origin 
through an
integral point. The duality between faces of $P$ and cones 
of $\Sigma$
is given as follows: For each face $E$ of $P$ let ${\Cal 
F}_E$ be the
set of codimension one faces containing $E$, and let
$$
\check\sigma _E =\SPAN \{ n _F \mid F \in {\Cal F}_E\} \ ,
$$
where $n _F$ is a primitive integral point, orthogonal to 
$F$,
and $n _F \cdot (m-p)>0$ for $p\in F$ and $m$ a vertex not 
in $F$.
(For $E=P$, $\breve \sigma _E =\{0\}$.) For example, 
$CP^n$ is the
toric variety determined by the simplex spanned by the unit
orthogonal basis vectors $e _1,\ldots, e_n$
and the origin or by the complete
fan whose rays are the half lines spanned by $e _1,\ldots, 
e_n$,
$-e _1 -e_2 -\cdots - e_n$.

For simplicity we assume our toric variety $X$ is 
simplicial;
i.e., each real $n$-dimensional cone in $\Sigma$ is 
spanned by $n$
rays. The present results can be modified to cover the 
general
case. Such an $X$ has rational singularities and can be 
described
locally as the quotient of ${\Bbb C}^n$ by the linear 
action of a 
finite cyclic group. In particular, $X$ is an orbifold and 
satisfies
Poincar\'e duality over ${\Bbb Q}$. A toric variety is 
smooth
if and only if it is simplicial and the integral points on
the rays of each cone span over ${\Bbb Z}$ a summand of 
the lattice
of integral points of ${\Bbb R}^n$.

For $X$ simplicial we give an explicit formula for $Td(X)
= T _0(X) =$\break $IT_0(X)\in H_\ast (X;{\Bbb C})$. With 
Riemann-Roch (see [D]
or the exposition [F2]) this leads to a formula for the 
number
of integral points in $P$. Actually, we will consider
$$
T(X) = [X] + 2 Td _{2n-2} (X) + \cdots + 2^n Td _0 (X)\ .
$$

Mock characteristic classes of $X$ may be defined as 
follows:
For each cone $\sigma\in \Sigma$, let $V(\sigma)$ be the 
corresponding
subvariety of $X$. If $\sigma =\breve\sigma _E$, $V(\sigma)$
is just the part of $X$ lying over the face $E$ and is, in 
fact,
the corresponding simplicial toric variety.
Now pretend the stable complex tangent bundle is a sum of 
line bundles 
with Chern classes $[V(\sigma)]\in H_{2n-2}(X;{\Bbb 
Q})\cong H^2(X;{\Bbb Q})$,
for $\sigma \in\Sigma^{(1)}$ a ray, and take the 
appropriate polynomial in 
these, as in [H, \S12] for smooth varieties. For example, 
for the
Todd polynomial this gives the Mock Todd class $TD(X)$ of 
[P]
(thought of in homology). Similarly, the mock $L$-classes
$$
L^{(m)} (X)\in H _\ast (X;{\Bbb Q})
$$
are defined using the $L$-polynomial. For smooth toric 
varieties,
the mock classes agree with the actual ones. For example, 
for $CP^n$
this gives the well-known formulas
$$
L(CP^n)=(c^2/\tanh c^2)^{n+1}
$$
and
$$
Td(CP^n)=\{ c/(1-e^{-c})\}^{n+1}\ ,
$$
$c$ represented by a hyperplane $CP^{n-1}\subset CP^n$.
Our mock $T$-class will be defined as
a sum over all cones $\sigma$,
$$
T^{(m)}(X)=\sum _\sigma L^{(m)}(V(\sigma))\ .
$$

Next we define some algebraic numbers associated to a 
$k$-dimensional
simplicial cone $\sigma$, essentially measuring ``angles'' 
in which
faces meet. Let $N\subset {\Bbb R}^n$ be the set of 
integral points.
By replacing $N$ by the plane containing $\sigma$ (or 
dividing by
its orthogonal complement), it may be assumed without loss 
of
generality that $n=\dim \, N = k$. Let $\sigma =\langle 
n_1,\ldots,
n_k\rangle$ be spanned by $n _i\in N$, $1\le i\le k$. Let 
$m _i$
be the unique primitive elements of $N$ with
$$
m _i \cdot n _j =0\ , \qquad i \ne j \ ,
$$
and
$$
q _i = m _i \cdot n _ i > 0 \ .
$$
Let $G _\sigma$ be the finite abelian group
$$
G _\sigma = N / N' \ , \qquad N' = {\Bbb Z} m_1 +\cdots +
{\Bbb Z} m_k \ .
$$
Set $m(\sigma)=\text{order of } G_\sigma$.
These are the groups which  act on affine spaces in the 
orbifold
description of X.

For $g =m + N'$ define $\lambda _j (g)$, also denoted 
$\lambda _{n_j}(g)$,
$$
\lambda _j (g) =\exp {2\pi i\over q_j} m 
\cdot m _j =\exp 2\pi i
\, \gamma _{n_j}(g)\ .
$$
Let
$$
G^o_\sigma =\{ g \in G _\sigma\mid \lambda _j (g)\ne 1
\text{ for } 1\le j\le k\} \ .
$$
Equivalently, $G^o_\sigma$ consists of the elements of $G 
_\sigma$
of the form $m +N'$ with $m$ in the interior of the cone 
spanned
by $m_1,\ldots, m_k$.

Now suppose that $\sigma$ as above is a cone of the 
simplicial
fan $\Sigma$. Let $\rho _i$ be the ray through $n _i$. 
Then we define
a characteristic class ${\Cal A}(\sigma)\in H_\ast 
(X_\Sigma;{\Bbb Q})$ by
$$
{\Cal A}(\sigma)={1\over m(\sigma)} \sum_{G^o_\sigma}
\prod ^k _1 {\lambda _i (g) e^{2[V(\rho_i)]}+1\over 
\lambda _i (g)
e^{2[V(\rho_i)]}-1} \ .
$$
If $G^o_\sigma$ is empty, we take ${\Cal A}(\sigma)=0$, 
and we also set
${\Cal A}(\sigma)=1$ for $\dim\, \sigma = 0$. The fact 
that this is
a rational class follows from Galois invariance.

\proclaim{Theorem 4}
Let $X _\Sigma$ be the toric variety corresponding to the 
complete
simplicial fan $\Sigma$. Then
$$
T(X _\Sigma) =\sum _\sigma T^{(m)}(V(\sigma))\cdot {\Cal 
A}(\sigma)\ .
$$
\endproclaim

This theorem applied in dimension $2n-4$ gives precisely 
Pommersheim's
calculation [P] of $Td_{2n-4}$ (denoted $Td^2$ in [P]).
Note that the present methods are quite different from 
those of [P].

Theorem~4 really is an explicit calculation, because the 
ring
$H_\ast (X;{\Bbb Q})$
$\cong$
$H^\ast (X;{\Bbb Q})$
$\cong 
A^\ast (X)\otimes {\Bbb Q}$ is explicitly known (see [O, 
pp.~134--135]).
Therefore, this result allows us to determine the Hilbert 
polynomial
(also called the Ehrhart polynomial) $\ell _P$, a degree $n$
polynomial [E] that counts lattice points in integral 
dilations
of $P$:
$$
\ell _P(\mu)=\bigsharp \{ N\cap \mu P\} \ , \qquad
\mu\text{ a positive integer.}
$$
Here is the result for $P$ a simplex, extending Pick's 
theorem
\cite{Pi} on
planar polygons and results of Mordell [M] and Pommersheim 
[P] in
dimension three to all dimensions: notation like $m(E)$, 
etc.,
refers to $m(\breve \sigma _E)$; $\nu (E)$ is the volume 
of the
face $E$, normalized with respect to the intersection of 
$N$ 
with $\dim \, E$ plane containing $E$; and 
${\Cal H}_E ={\Cal F}_\phi-{\Cal F}_E$.

\proclaim{Theorem 5}
Let $\Delta$ be an $n$-simplex with vertex points in the 
lattice
$N$. Let $a_ r$ be the coefficient of $U^r$ in the power 
series
$$
\sum_{E\le\Delta }{1\over m (E)} \left[ \sum_{K\le E}
\omega _K \prod_{F\in {\Cal H}_K} {(\nu(F)U)\over 
\tanh(\nu(F)U)}\right]
\sum_{G^o_E} \prod_{F\in {\Cal F}_E} \coth\left\{
\pi i \gamma^E _F(g) +\nu(F)U\right\} \ ,
$$
where
$$
\omega _K = m(K)\prod_{F\in {\Cal F}_K}(\nu(F)U) \ .
$$
Let $E$ be any face of $\Delta$ of dimension $r$, and let
$$
b _r = {\nu(E)\over 2^{n-r}m(E)\prod_{F\in {\Cal 
F}_E}\nu(F)} \ .
$$
Then
$$
\ell _\Delta (k)=\sum^n_{r=0} a_r b_r k^r\ .
$$
\endproclaim

These results and known facts about Todd classes or Ehrhart
polynomials imply nontrivial number theoretical results
(cf.\ \cite{R1, R2, P}). For other valuable perspectives
and results on
Todd classes of toric varieties and lattice points, see the
works of
Brion [B] and of Morelli [Mo].

Note that in this formula all terms are just functions of 
the volumes
of the codimension one faces, the volume of a single face in
each dimension, and the cones at each vertex.

The ingredients of the proof of Theorem~4 are Theorem~2 
applied
to a resolution of $X$ obtained by interior subdivision of
the fan $\Sigma$, formulas of [H] relating Todd classes and
$L$-classes of smooth varieties, and formulas of 
Atiyah-Singer
and Hirzebruch-Zagier type for equivariant $L$-classes 
related
to the orbifold structure.

\noindent{\bf Added May 17, 1993.}
The first two coefficients of $\ell _\Delta$ (and the 
constant term)
are classically known (essentially ``Pick's Theorem''),
and Pommersheim's work mentioned above determines the next
coefficient. Khovanskii has just communicated to us
that recently but independently of Pommersheim (and of
the above results determining all the coefficients of
$\ell _\Delta)$ he and Kantor have also found one 
coefficient.
\frenchspacing



\enddocument